\documentclass[12pt,reqno]{article}

\usepackage[usenames]{color}
\usepackage{amssymb}
\usepackage{graphicx}
\usepackage{amscd}
\usepackage{lscape}

\usepackage{amsthm}
\newtheorem{theorem}{Theorem}

\theoremstyle{definition}

\newtheorem{example}[theorem]{Example}

\usepackage[colorlinks=true,
linkcolor=webgreen, filecolor=webbrown,
citecolor=webgreen]{hyperref}

\definecolor{webgreen}{rgb}{0,.5,0}
\definecolor{webbrown}{rgb}{.6,0,0}

\usepackage{color}

\usepackage{float}

\usepackage{graphics,amsmath,amssymb}
\usepackage{amsfonts}
\usepackage{latexsym}
\usepackage{epsf}

\newcommand{\seqnum}[1]{\href{http://www.research.att.com/cgi-bin/access.cgi/as/~njas/sequences/eisA.cgi?Anum=#1}{\underline{#1}}}

\begin{document}

\begin{center}
\vskip 1cm{\LARGE\bf A Note on $d$-Hankel Transforms, Continued Fractions, and Riordan Arrays} \vskip 1cm \large
Paul Barry\\
School of Science\\
Waterford Institute of Technology\\
Ireland\\
\href{mailto:pbarry@wit.ie}{\tt pbarry@wit.ie}
\end{center}
\vskip .2 in

\begin{abstract} The Hankel transform of an integer sequence is a much studied and much applied mathematical operation. In this note, we extend the notion in a natural way to sequences of $d$ integer sequences. We explore links to generalized continued fractions in the context of $d$-orthogonal sequences.\end{abstract}

\section{Introduction}
The Hankel transform of sequences \cite{Kratt, KrattC, Layman} is a topic that has attracted some attention in recent years. Many Hankel transforms have been shown to have combinatorial interpretations, and sequences with interesting Hankel transforms are often interesting in themselves. When the elements of the Hankel transform of a sequence are all positive, there is a classical link to orthogonal polynomials \cite{Chihara, Gautschi, Szego}.
Many Hankel transforms are documented in the On-Line Encyclopedia of Integer Sequences \cite{SL1, SL2}.

In this note, we look at the notion of the $d$-Hankel transform of a family of $d$ sequences. We use the theory of Riordan arrays \cite{Meixner, Hankel_Riordan, SGWW, Survey} to narrow the focus of the note to families of sequences that are defined by the columns of appropriate Riordan arrays. In this context, ``appropriate'' will mean that the Riordan arrays in question have production matrices \cite{ProdMat_0, ProdMat} that are $d$-diagonal.

The $d$-Hankel transforms in this note will be closely linked to the theory of $d$-orthogonal polynomials \cite{Lamiri, Maroni, Isegham1, Isegham2}. We recall that a family of polynomials $P_n(x)$, with $P_n(x)$ of precise degree $n$, is a family of $d$-orthogonal polynomials \cite{Lamiri, Maroni} if $P_{n+1}$ can be specified by the $d+1$ recurrence
$$xP_n(x)=\sum_{k=0}^{d+1} \alpha_{k,d}(n)P_{n-d+k}(x),$$
where $\alpha_{d+1,k} \alpha_{0,d}(n) \ne 0$, $n\ge 0$, where by convention, $P_{-n}=0, n\ge 1$.

In order to motivate our discussion, we begin with an example. For this, we use the language of exponential Riordan arrays, wherein an exponential Riordan array is defined by two power series
$$g(x)=1+g_1 \frac{x}{1!}+g_2 \frac{x^2}{2!}+ \cdots = \sum_{n=0}^{\infty}g_n \frac{x^n}{n!},$$ and
$$f(x)=f_1 \frac{x}{1!}+f_2 \frac{x^2}{2!}+\cdots =\sum_{n=1}^{\infty} f_n \frac{x^n}{n!}.$$ We require that $f_1 \ne 0$ and we often assume that $f_1=1$.
The $(n,k)$-term of the lower-triangular matrix defined by the pair $(g(x), f(x))$ is then given by
$$T_{n,k}=\frac{n!}{k!} [x^n] g(x) f(x)^k.$$ The operator $[x^n]$ here extracts the coefficient of $x^n$ \cite{MC}.

This infinite lower-triangular matrix is denoted by $[g(x), f(x)]$, and it is called the exponential Riordan array defined by the pair of power series $(g(x), f(x))$. Such a matrix is invertible, and we have
$$[g(x), f(x)]^{-1} = \left[\frac{1}{g(\bar{f}(x))}, \bar{f}(x)\right],$$ where
$\bar{f}(x)$ is the compositional inverse of $f(x)$. This means that we have
$$\bar{f}(f(x))=x,\quad\quad f(\bar{f}(x))=x.$$
If $A=[g(x), f(x)]$ is an exponential Riordan array, an important object related to $A$ is its production matrix $P=P_A$, defined by
$$P = A^{-1} \bar{A},$$ where $\bar{A}$ is the matrix $A$ with its top line removed.
It can be shown that the matrix $P_A$ has as bivariate generating function the expression
$$ e^{xy} (Z(x) + y A(x)),$$ where
$$A(x)=f'(\bar{f}(x))=\frac{1}{\bar{f}'(x)}$$ and
$$Z(x)=\frac{g'(\bar{f}(x))}{g(\bar{f}(x))}.$$

\begin{example}
We consider the exponential Riordan array
$$\left[\frac{1}{\sqrt{1-2x}}, \frac{1}{\sqrt{1-2x}-1}\right].$$
This array begins
$$\left(
\begin{array}{ccccccc}
 1 & 0 & 0 & 0 & 0 & 0 & 0 \\
 1 & 1 & 0 & 0 & 0 & 0 & 0 \\
 3 & 5 & 1 & 0 & 0 & 0 & 0 \\
 15 & 33 & 12 & 1 & 0 & 0 & 0 \\
 105 & 279 & 141 & 22 & 1 & 0 & 0 \\
 945 & 2895 & 1830 & 405 & 35 & 1 & 0 \\
 10395 & 35685 & 26685 & 7500 & 930 & 51 & 1 \\
\end{array}
\right).$$
Its production matrix, which has $Z(x)=(1+x)^2$ and $A(x)=(1+x)^3$,  begins
$$\left(
\begin{array}{ccccccc}
 1 & 1 & 0 & 0 & 0 & 0 & 0 \\
 2 & 4 & 1 & 0 & 0 & 0 & 0 \\
 2 & 10 & 7 & 1 & 0 & 0 & 0 \\
 0 & 12 & 24 & 10 & 1 & 0 & 0 \\
 0 & 0 & 36 & 44 & 13 & 1 & 0 \\
 0 & 0 & 0 & 80 & 70 & 16 & 1 \\
 0 & 0 & 0 & 0 & 150 & 102 & 19 \\
\end{array}
\right).$$
The $4$-diagonal nature of the production matrix indicates that the inverse matrix of $A$, given by
$$\left[\frac{1}{\sqrt{1-2x}}, \frac{1}{\sqrt{1-2x}-1}\right]^{-1}=\left[\frac{1}{1+x}, \frac{1}{2}\left(1-\frac{1}{(1+x)^2}\right)\right],$$ is the coefficient array of the family of polynomials $P_n(x)$ which begins
$$1, x-1, x^2-5x+2, x^3-12x^2+27x-6, x^4-22x^3+123x^2-168x+24,\ldots$$ and which constitutes a family of
$2$-orthogonal polynomials. This means that they obey a recurrence of the form
$$P_n(x)=(x-\alpha_n) P_{n-1}(x)-\beta_n P_{n-2}(x)-\gamma_n P_{n-3}(x).$$ In this case, the values of
$\alpha_n, \beta_n$ and $\gamma_n$ may be read from the diagonal elements of the production matrix.

We now look at the first two columns of $A$, given respectively by
$$1, 1, 3, 15, 105, 945, 10395, 135135,\ldots \quad\text{\seqnum{A001147}}$$ with e.g.f. given by $$\frac{1}{\sqrt{1-2x}},$$ and
$$0, 1, 5, 33, 279, 2895, 35685, 509985, \ldots$$ with e.g.f. given by
$$\frac{1}{\sqrt{1-2x}}\left(\frac{1}{\sqrt{1-2x}}-1\right).$$
We now form the sequence that is the sum of these two sequences, giving us the sequence that begins
$$1, 2, 8, 48, 384, 3840, 46080, 645120,\ldots \quad\text{\seqnum{A0000165}}$$ with e.g.f. given by
$$\frac{1}{\sqrt{1-2x}}+\frac{1}{\sqrt{1-2x}}\left(\frac{1}{\sqrt{1-2x}}-1\right)=\frac{1}{1-2x}.$$ This latter sequence is in fact $2^n n!$, while the first sequence is $(2n-1)!!$.
We let $a_n = (2n-1)!!$ and $b_n=2^n n!$. We define the 2-Hankel transform of $(a_n, b_n)$ to be the sequence
$h_n$ defined by the determinant values
$$h_n = \left| \begin{cases}
   a_{i+j-\lfloor \frac{i}{2} \rfloor} & \text{if $i$ is even} \\
   b_{i+j-\lfloor \frac{i+1}{2} \rfloor} & \text{if $i$ is odd}
   \end{cases} \right|_{0 \le i,j \le n}$$

We find that the $2$-Hankel transform of $((2n-1)!!, 2^n n!)$ begins
$$1, 1, 2, 24, 1728, 1658880, 17915904000, 4334215495680000, \ldots$$
which is equal to the product
$$\prod_{k=0}^n \gamma_k^{\lfloor \frac{n-k}{2} \rfloor}=\prod_{k=0}^n ((k+2)(k+1)^2)^{\lfloor \frac{n-k}{2} \rfloor},$$ where
$\gamma_n$ is the sequence $(n+2)(n+1)^2$ that begins
$$2, 12, 36, 80, 150, 252, 392, 576,\ldots.$$
We note that if we replace $a_n$ and $b_n$ by their binomial transforms
$\sum_{k=0}^n \binom{n}{k}a_k$ and $\sum_{k=0}^n \binom{n}{k}b_k$, respectively, we obtain the same
$2$-Hankel transform.

We can recover the polynomials $P_n(x)$ using determinants as follows. We have

$$P_n(x)= \frac{h_n(x)}{h_{n-1}}$$ where
$h_n(x)$ is the same as the determinant $h_n$, except that the last row is given by
$1,x,x^2, \ldots$.

\end{example}
\begin{example}
Our next example is based on the ordinary Riordan array $A$ with general term given by
$$T_{n,k}= \binom{4n} {3n+k}.$$ We recall that an ordinary Riordan array, or Riordan array for short, is defined by two power series
$$g(x)=1+g_1 x+ g_2 x +\cdots = \sum_{n=0}^{\infty} g_n x^n$$ and
$$f(x)= x + f_2 x^2+ f_3x^3 + \cdots = \sum_{n=1}^n f_n x^n.$$ Then
we have
$$T_{n,k}=[x^n] g(x) f(x)^k.$$
In this case, the array $A$ begins
$$\left(
\begin{array}{ccccccc}
 1 & 0 & 0 & 0 & 0 & 0 & 0 \\
 4 & 1 & 0 & 0 & 0 & 0 & 0 \\
 28 & 8 & 1 & 0 & 0 & 0 & 0 \\
 220 & 66 & 12 & 1 & 0 & 0 & 0 \\
 1820 & 560 & 120 & 16 & 1 & 0 & 0 \\
 15504 & 4845 & 1140 & 190 & 20 & 1 & 0 \\
 134596 & 42504 & 10626 & 2024 & 276 & 24 & 1 \\
\end{array}
\right),$$ and it has a production matrix $A^{-1}\bar{A}$ that begins
$$\left(
\begin{array}{ccccccc}
 4 & 1 & 0 & 0 & 0 & 0 & 0 \\
 12 & 4 & 1 & 0 & 0 & 0 & 0 \\
 12 & 6 & 4 & 1 & 0 & 0 & 0 \\
 4 & 4 & 6 & 4 & 1 & 0 & 0 \\
 0 & 1 & 4 & 6 & 4 & 1 & 0 \\
 0 & 0 & 1 & 4 & 6 & 4 & 1 \\
 0 & 0 & 0 & 1 & 4 & 6 & 4 \\
\end{array}
\right).$$
This means that the elements of $A^{-1}$ are the coefficients of a family $P_n(x)$ of $3$-orthogonal polynomials, that begins
$$1, x-4, x^2-8x+4, x^3-12x^2+30x-4, x^4-16x^3+72x^2-80x +4,\ldots.$$
We have
$$P_n(x)=(x-4)P_{n-1}(x)-6 P_{n-2}(x)-4 P_{n-3}(x)-P_{n-4}(x),$$ where
$P_0(x)=1$, $P_1(x)=x-4$, $P_2(x)=x^2-8x+4$, and $P_3(x)=x^3-12x^2+30x-4$.

We now take the first column of $A$, the sum of the first two columns of $A$, and the sum of the first three columns of $A$ to obtain three sequences $a_n$, $b_n$ and $c_n$ that begin, respectively, as follows.
$$a_n : 1, 4, 28, 220, 1820, 15504, 134596, 1184040,\ldots$$
$$b_n : 1, 5, 36, 286, 2380, 20349, 177100, 1560780,\ldots$$
$$c_n : 1, 5, 37, 298, 2500, 21489, 187726, 1659060,\ldots$$
We then define the $3$-Hankel transform $h_n$ of these sequences as follows.

$$h_n = \left| \begin{cases}
   a_{i+j-\lfloor \frac{2i}{3} \rfloor} &  i \bmod 3 = 0 \\
   b_{i+j-\lfloor \frac{2i+1}{3} \rfloor} & i \bmod 3 = 1\\
   c_{i+j-\lfloor \frac{2i+2}{3} \rfloor} & i \bmod 3 =2 \\
   \end{cases} \right|_{0 \le i,j \le n}$$
We find that $h_n$ begins
$$1, 1, 1, 4, 4, 4, 16, 16, 16,\ldots$$ corresponding to the sequence $\gamma_n$ that begins
$$4,1,1,1,1,1,\ldots.$$
\end{example}
\section{Continued fractions}
We now look at the continued fraction for the sequence $(2n-1)!!$. It is customary to associate continued fractions of Jacobi or Stieltjes type \cite{Wall} with the moments of orthogonal polynomials. The form of generalized continued fraction that we need in this note has been studied, for instance, in the context of lattice path enumeration \cite{Varvak}.  The production matrix of the last section infers that the generating function of $(2n-1)!!$  has the following form.
$$\cfrac{1}{1-x-\cfrac{2x^2}{1-\cdots}-\cfrac{2x^3}{(1-\cdots)(1-\cdots)}}.$$
See the Appendix for a more explicit form.

Similarly, the continued fraction for $\binom{4n}{3n}$ takes the following form.
$$\cfrac{1}{1-4x-\cfrac{12x^2}{1-\cdots}-\cfrac{12x^3}{(1-\cdots)(1-\cdots)}-\cfrac{4x^4}{(1-\cdots)(1-\cdots)(1-\cdots)}}.$$

\section{Constructing $d$ orthogonal polynomials from Riordan arrays}
For an ordinary Riordan array to represent the coefficient array of a family of $d$ orthogonal polynomials,
we require that the production array of its inverse have a $Z$ sequence and an $A$ sequence of the type
$$Z(x)=\sum_{i=0}^d Z_i x^i,$$ and
$$A(x)=\sum_{i=0}^{d+1} A_i x^i.$$
The inverse array (with $Z$-sequence $Z(x)$ and $A$-sequence $A(x)$)  is given by
$$\left(\frac{1}{1-x Z\left(\text{Rev}\left\{\frac{x}{A(x)}\right\}\right)}, \text{Rev}\left\{\frac{x}{A(x)}\right\}\right),$$
while the coefficient matrix of the family of $d$ orthogonal polynomials will be given by
$$\left(1-\frac{xZ}{A}, \frac{x}{A}\right).$$
We recall that in the case of an ordinary Riordan array $(g(x), f(x))$, we have
$$Z(x)= \frac{1}{\bar{f}(x)}\left(1-\frac{1}{g(\bar{f}(x))}\right),$$ and
$$A(x)=\frac{x}{\bar{f}(x)}.$$
\begin{example} We consider the Riordan array $(g(x), f(x))$ with
$$Z(x)=1+x+x^2$$ and
$$A(x)=1+x+2x^2+3x^3.$$
Then the coefficient array of the corresponding $2$-orthogonal family is given by
$$\left(1-\frac{xZ}{A}, \frac{x}{A}\right)=\left(1-\frac{x(1+x+x^2)}{1+x+2x^2+3x^3}, \frac{x}{1+x+2x^2+3x^2}\right),$$ or $$\left(\frac{1+x^2+2x^3}{1+x+2x^2+3x^3}, \frac{x}{1+x+2x^2+3x^3}\right).$$
This matrix begins
$$\left(
\begin{array}{ccccccc}
 1 & 0 & 0 & 0 & 0 & 0 & 0 \\
 -1 & 1 & 0 & 0 & 0 & 0 & 0 \\
 0 & -2 & 1 & 0 & 0 & 0 & 0 \\
 1 & 0 & -3 & 1 & 0 & 0 & 0 \\
 2 & 2 & 1 & -4 & 1 & 0 & 0 \\
 -4 & 6 & 4 & 3 & -5 & 1 & 0 \\
 -3 & -14 & 9 & 6 & 6 & -6 & 1 \\
\end{array}
\right)$$ and hence gives rise to the $2$-orthogonal family of polynomials $P_n(x)$ that begins
$$1, x-1, x^2-2x, x^3-3x^2+1, x^4-4x^3+x^2+2, x^5 - 5x^4 + 3x^3 + 4x^2 + 6x - 4,\ldots$$ and that obeys the following
$4$-term recurrence.
$$P_n(x)=(x-1)P_{n-1}(x)-2 P_{n-2}(x)-3 P_{n-3}(x),$$ with $P_0(x)=1$, $P_1(x)=x-1$ and $P_2(x)=x^2-2x$.

The inverse of the coefficient array begins
$$\left(
\begin{array}{ccccccc}
 1 & 0 & 0 & 0 & 0 & 0 & 0 \\
 1 & 1 & 0 & 0 & 0 & 0 & 0 \\
 2 & 2 & 1 & 0 & 0 & 0 & 0 \\
 5 & 6 & 3 & 1 & 0 & 0 & 0 \\
 14 & 20 & 11 & 4 & 1 & 0 & 0 \\
 45 & 68 & 42 & 17 & 5 & 1 & 0 \\
 155 & 248 & 159 & 72 & 24 & 6 & 1 \\
\end{array}
\right),$$ which therefore has production matrix
$$\left(
\begin{array}{ccccccc}
 1 & 1 & 0 & 0 & 0 & 0 & 0 \\
 1 & 1 & 1 & 0 & 0 & 0 & 0 \\
 1 & 2 & 1 & 1 & 0 & 0 & 0 \\
 0 & 3 & 2 & 1 & 1 & 0 & 0 \\
 0 & 0 & 3 & 2 & 1 & 1 & 0 \\
 0 & 0 & 0 & 3 & 2 & 1 & 1 \\
 0 & 0 & 0 & 0 & 3 & 2 & 1 \\
\end{array}
\right).$$
We thus find that the $2$-Hankel transform of the sequences $a_n$
$$1, 1, 2, 5, 14, 45, 155, 562, 2122, 8245, 32769,\ldots$$
and $b_n$
$$1, 2, 4, 11, 34, 113, 403, 1499, 5758, 22691, 91189,\ldots$$
is given by
$$\prod_{k=0}^n (3-2*0^k)^{\lfloor \frac{n-k}{2} \rfloor},$$ which begins
$$1, 1, 1, 3, 9, 81, 729, 19683, 531441, 43046721, 3486784401,\ldots.$$

Note that in this case, the ``moment'' matrix $(g(x), f(x))$ has
$$f(x)=\frac{2}{9}\left(\sqrt{\frac{9-5x}{x}}\sin\left(\frac{1}{3}\sin^{-1}\left(\frac{x(205x+54)}{2(5x-9)^2}\sqrt{\frac{9-5x}{x}}\right)\right)-1\right).$$
The corresponding $g(x)$ is given in the Appendix. The sequence $a_n$ has g.f. $g(x)$ while the sequence $b_n$ has g.f. $g(x)(1+f(x))$. (We call the inverse of the coefficient array of a family of $d$-orthogonal polynomials a ``moment'' matrix, in analogy to the situation for orthogonal polynomials).
\end{example}
\begin{example} We now turn to an example of a family of $2$-orthogonal polynomials defined by an exponential Riordan array. The coefficient array of this family will therefore be the inverse of an exponential Riordan array
$$[g(x), f(x)].$$ The matrix $[g(x), f(x)]$ must therefore have a $Z$ sequence and an $A$-sequence of the appropriate form. We recall that we have
$$[g(x), f(x)] = \left[e^{\int_0^{\text{Rev}\left(\int_0^x \frac{dt}{A(t)}\right)}\frac{Z(t)}{A(t)}\,dt}, \text{Rev}\left(\int_0^x \frac{dt}{A(t)}\right)\right],$$ while
$$[g(x), f(x)]^{-1} = \left[\frac{1}{e^{\int_0^x \frac{Z(t)}{A(t)}\,dt}}, \int_0^x \frac{dt}{A(t)}\right].$$
We revisit our first example in light of these facts.
Thus we assume that
$$A(x)=1+3x+3x^2+x^3=(1+x)^3,$$ and that
$$Z(x)=1+2x+x^2=(1+x)^2.$$
We recall that the production matrix $P$ will now have generating function (exponential in $x$, ordinary in $y$) given by
$$e^{xy}((1+x)^2+y(1+x)^3).$$
We have
$$\int_0^x \frac{dt}{(1+t)^3} = \frac{1}{2}\left(1-\frac{1}{(1+x)^2}\right)=\frac{x(2+x)}{2(1+x)^2},$$
and
$$\int_0^x \frac{Z(t)}{A(t)}\,dt=\int_0^x \frac{(1+t)^2}{(1+t)^3}\,dt=\ln(1+x).$$
Hence we find that the coefficient array sought is given by
$$[g(x), f(x)]^{-1}=\left[ \frac{1}{1+x}, \frac{1}{2}\left(1-\frac{1}{(1+x)^2}\right)\right],$$ with therefore
$$[g(x), f(x)]=\left[\frac{1}{\sqrt{1-2x}}, \frac{1}{\sqrt{1-2x}}-1\right].$$
\end{example}
\begin{example}
In this example, we look at the following $3$-orthogonal situation.
We choose the exponential Riordan array $[g(x), f(x)]$ such that
$$A(x)=1+x+x^2+x^3+x^4$$ and
$$Z(x)=1+2x+3x^2+4x^3.$$ Then we have
$$\int_0^x \frac{1}{A(t)}=\bar{f}(x) \quad\text{(see the Appendix)},$$ while
$$\int_0^x \frac{Z(t)}{A(t)}\,dt = \ln(1+x+x^2+x^3+x^4),$$ and hence we have
$$[g(x), f(x)]^{-1}=\left[\frac{1}{1+x+x^2+x^3+x^4}, \bar{f}(x)\right].$$
We can constitute the triangle $[g(x), f(x)]$ by using the production matrix whose generating function is given by
$$e^{xy}(1+2x+3x^2+4x^3+y(1+x+x^2+x^3+x^4)).$$ We find that the production matrix begins
$$\left(
\begin{array}{ccccccc}
 1 & 1 & 0 & 0 & 0 & 0 & 0 \\
 2 & 2 & 1 & 0 & 0 & 0 & 0 \\
 6 & 6 & 3 & 1 & 0 & 0 & 0 \\
 24 & 24 & 12 & 4 & 1 & 0 & 0 \\
 0 & 120 & 60 & 20 & 5 & 1 & 0 \\
 0 & 0 & 360 & 120 & 30 & 6 & 1 \\
 0 & 0 & 0 & 840 & 210 & 42 & 7 \\
\end{array}
\right).$$
This in turn generates the matrix $[g(x), f(x)]$ which begins
$$\left(
\begin{array}{ccccccc}
 1 & 0 & 0 & 0 & 0 & 0 & 0 \\
 1 & 1 & 0 & 0 & 0 & 0 & 0 \\
 3 & 3 & 1 & 0 & 0 & 0 & 0 \\
 15 & 15 & 6 & 1 & 0 & 0 & 0 \\
 105 & 105 & 45 & 10 & 1 & 0 & 0 \\
 825 & 945 & 420 & 105 & 15 & 1 & 0 \\
 7755 & 9555 & 4725 & 1260 & 210 & 21 & 1 \\
\end{array}
\right).$$
We now obtain the sequences formed by the first column, the sum of the first and second columns, and the sum of the first three columns, to get
$$a_n : 1, 1, 3, 15, 105, 825, 7755, 85455, 1076625, 15154425, \ldots$$
$$b_n : 1, 2, 6, 30, 210, 1770, 17310, 196110, 2531250, 36545850,\ldots$$
$$c_n : 1, 2, 7, 36, 255, 2190, 22035, 255120, 3351915, 49198050, \ldots$$
We find that the $3$-Hankel transform of these sequences, which begins
$$1, 1, 1, 24, 2880, 1036800, 20901888000, 4213820620800000, 4587333680627712000000,\ldots,$$ is equal to
$$\prod_{k=0}^n (24 \binom{k+4}{4})^{\lfloor \frac{n-k}{3} \rfloor}=\prod_{k=0}^n ((k+1)(k+2)(k+3)(k+4))^{\lfloor \frac{n-k}{3} \rfloor}.$$
\end{example}
\section{Recurrence coefficients}
As in the case of orthogonal polynomials, it is possible to find determinantal expressions for the recurrence coefficients. As $d$ increases, these formulas grow more complex. We restrict ourselves to looking at the case $d=2$. Thus we seek formulas for the  coefficients in the recurrence
$$P_n(x)=(x-\alpha_n)P_{n-1}(x)- \beta_n P_{n-2}(x)-P_{n-3}(x).$$
We know that we also have $$P_n(x)=\frac{h_n(x)}{h_{n-1}}.$$ Substituting similar expressions for $P_{n-1}(x)$, $P_{n-2}(x)$ and $P_{n-3}(x)$ into the $4$-term recurrence, and equating coefficients of powers of $x$, we obtain that
$$\alpha_n = \frac{h_{n+1,n}}{h_n}-\frac{h_{n,n-1}}{h_{n-1}}+0^n,$$
$$\beta_n = \frac{\alpha_n h_{n,n-1}}{h_{n-1}}+\frac{h_{n,n-2}}{h_{n-1}}-\frac{h_{n+1,n-1}}{h_n},$$
and
$$\gamma_n = \frac{\beta_n h_{n-1,n-2}}{h_{n-2}}-\frac{\alpha_n h_{n,n-2}}{h_{n-1}}-\frac{h_{n,n-3}}{h_{n-1}}+\frac{h_{n+1,n-2}}{h_n},$$ with appropriate adjustments for initial terms.
\section{Generalizations} In this note, we have confined our attention to the case $d=2$ and $d=3$. To define the $d$-Hankel transform, we can take the first $d$ column sums of a ``moment'' matrix that is generated by a suitable production matrix that is $d$-diagonal. We thus obtain a family $a_{n,k}$ of sequences, where $0 \le k \le d-1$. We can then define the $d$-Hankel transform of this family of $d$ sequences by
$$h_n^{(d)} = \left| \begin{cases}
   a_{i+j-\lfloor \frac{(d-1)i}{d} \rfloor,0} &  i \bmod d = 0 \\
   a_{i+j-\lfloor \frac{(d-1)i+1}{d} \rfloor,1} & i \bmod d = 1\\
   a_{i+j-\lfloor \frac{(d-1)i+2}{d} \rfloor,2} & i \bmod d =2 \\
   \ldots\\
   a_{i+j-\lfloor \frac{(d-1)i+d-1}{d} \rfloor,d-1} & i \bmod d = d-1
   \end{cases} \right|_{0 \le i,j \le n}$$
The corresponding family of $d$-orthogonal polynomials is then defined by
$$P_n(x)=\frac{h_n^{(d)}(x)}{h_{n-1}^{(d)}},$$
where $h_n^{(d)}(x)$ is the same as $h_n^{(d)}$ above, except that the last row is $1,x,x^2,\ldots$.

\section{Conclusions} In this short note, we have seen how the $d$-Hankel transform can be related to $d$-orthogonal polynomials, and we have used the particular context where the polynomials studied have Riordan arrays for coefficient arrays. Given the combinatorial interpretation of many Hankel transforms \cite{Jac, Bressoud}, it will be of interest to see which $d$-Hankel transforms also have strong combinatorial associations.

\pagebreak
\section{Appendix}
In this appendix we give a more explicit example of generalized continued fraction for the ``moments'' of a $2$-orthogonal polynomial family (given by $(2n-1)!!$). Also we give the explicit formulas for $g(x)$ and $\bar{f}(x)$ referenced in the examples above.

Given their size, these are to be found on the next page printed in landscape form.

\begin{landscape}
\begin{tiny}
$$\cfrac{1}{1-x-\cfrac{2x^2}{1-4x-\cfrac{10x^2}{1-7x-\cfrac{24x^2}{1-\cdots}-\cfrac{36x^3}{(1-\cdots)(1-\cdots)}}-\cfrac{12x^2}{(1-7x-\cfrac{24x^2}{1-\cdots}-\cfrac{36x^3}{(1-\cdots)(1-\cdots)})(1-10x-\cfrac{44x^2}{1-\cdots}-\cfrac{80x^3}{(1-\cdots)(1-\cdots)})}}-\cfrac{2x^3}{(1-4x-\cfrac{10x^2}{(\cdots)}-\cfrac{12x^2}{(\cdots)(\cdots)})(1-7x-\cfrac{24x^2}{(\cdots)}-\cfrac{36x^2}{(\cdots)(\cdots)})}}.$$
\end{tiny}
$$g(x)=\frac{54-54 \sqrt{\frac{9}{x}-5} \sin \left(\frac{1}{3} \sin ^{-1}\left(\frac{\sqrt{\frac{9}{x}-5} x (205
   x+54)}{2 (9-5 x)^2}\right)\right)}{-123 x+4 \sqrt{\frac{9}{x}-5} (7 x-9) \sin \left(\frac{1}{3} \sin
   ^{-1}\left(\frac{\sqrt{\frac{9}{x}-5} x (205 x+54)}{2 (9-5 x)^2}\right)\right)+2 (5 x-9) \cos
   \left(\frac{2}{3} \sin ^{-1}\left(\frac{\sqrt{\frac{9}{x}-5} x (205 x+54)}{2 (9-5 x)^2}\right)\right)+54}.$$
\begin{scriptsize}
$$\bar{f}(x)=-\frac{\log \left(\frac{2 x^2+\left(1-\sqrt{5}\right) x+2}{2 x^2+\left(1+\sqrt{5}\right) x+2}\right)}{2
   \sqrt{5}}+\sqrt{\frac{1}{5}-\frac{2}{5 \sqrt{5}}} \tan ^{-1}\left(\sqrt{\frac{1}{8}-\frac{1}{8 \sqrt{5}}}
   \left(4 x-\sqrt{5}+1\right)\right)+\sqrt{\frac{1}{5}+\frac{2}{5 \sqrt{5}}} \tan
   ^{-1}\left(\sqrt{\frac{1}{8}+\frac{1}{8 \sqrt{5}}} \left(4
   x+\sqrt{5}+1\right)\right)-\sqrt{\frac{1}{50}+\frac{1}{50 \sqrt{5}}} \pi.$$
   \end{scriptsize}
\end{landscape}

\bigskip
\hrule
\bigskip
\noindent 2010 {\it Mathematics Subject Classification}: Primary
15B36; Secondary 33C45, 11B83, 11C20, 05A15.
\noindent \emph{Keywords:} Riordan array, orthogonal polynomials, $d$-Hankel transform, $d$-orthogonal polynomials, generating functions.

\end{document}